\documentclass[12pt,reqno]{amsart}

\newtheorem{theorem}{Theorem}[section]
\newtheorem{lemma}[theorem]{Lemma}

\usepackage{amsmath,amssymb,amsfonts,}
\usepackage[dvips]{graphics}
\usepackage{amssymb}
\usepackage{amsfonts}
\usepackage{latexsym}
\usepackage{cite}
\usepackage{mathrsfs}

\theoremstyle{definition}

\theoremstyle{remark}
\newtheorem{remark}[theorem]{Remark}

\theoremstyle{conjecture}
\newtheorem{conjecture}[theorem]{Conjecture}

\theoremstyle{corollary}
\newtheorem{corollary}[theorem]{Corollary}

\theoremstyle{proposition}

\theoremstyle{problem}

\theoremstyle{plain}
\newtheorem{claim}{Claim}

\numberwithin{equation}{section}

\large\normalsize
\RequirePackage[normalem]{ulem} 
\RequirePackage{color}\definecolor{RED}{rgb}{1,0,0}\definecolor{BLUE}{rgb}{0,0,1} 
\providecommand{\DIFaddbegin}{} 
\providecommand{\DIFaddend}{} 
\providecommand{\DIFdelbegin}{} 
\providecommand{\DIFdelend}{} 
\providecommand{\DIFaddbeginFL}{} 
\providecommand{\DIFaddendFL}{} 
\providecommand{\DIFdelbeginFL}{} 
\providecommand{\DIFdelendFL}{} 
\newcommand{\DIFscaledelfig}{0.5}
\RequirePackage{settobox} 
\RequirePackage{letltxmacro} 
\newsavebox{\DIFdelgraphicsbox} 
\newlength{\DIFdelgraphicswidth} 
\newlength{\DIFdelgraphicsheight} 
\LetLtxMacro{\DIFOincludegraphics}{\includegraphics} 
\newcommand{\DIFaddincludegraphics}[2][]{{\color{blue}\fbox{\DIFOincludegraphics[#1]{#2}}}} 
\newcommand{\DIFdelincludegraphics}[2][]{
\sbox{\DIFdelgraphicsbox}{\DIFOincludegraphics[#1]{#2}}
\settoboxwidth{\DIFdelgraphicswidth}{\DIFdelgraphicsbox} 
\settoboxtotalheight{\DIFdelgraphicsheight}{\DIFdelgraphicsbox} 
\scalebox{\DIFscaledelfig}{
\parbox[b]{\DIFdelgraphicswidth}{\usebox{\DIFdelgraphicsbox}\\[-\baselineskip] \rule{\DIFdelgraphicswidth}{0em}}\llap{\resizebox{\DIFdelgraphicswidth}{\DIFdelgraphicsheight}{
\setlength{\unitlength}{\DIFdelgraphicswidth}
\begin{picture}(1,1)
\thicklines\linethickness{2pt} 
{\color[rgb]{1,0,0}\put(0,0){\framebox(1,1){}}}
{\color[rgb]{1,0,0}\put(0,0){\line( 1,1){1}}}
{\color[rgb]{1,0,0}\put(0,1){\line(1,-1){1}}}
\end{picture}
}\hspace*{3pt}}} 
} 
\LetLtxMacro{\DIFOaddbegin}{\DIFaddbegin} 
\LetLtxMacro{\DIFOaddend}{\DIFaddend} 
\LetLtxMacro{\DIFOdelbegin}{\DIFdelbegin} 
\LetLtxMacro{\DIFOdelend}{\DIFdelend} 
\DeclareRobustCommand{\DIFaddbegin}{\DIFOaddbegin \let\includegraphics\DIFaddincludegraphics} 
\DeclareRobustCommand{\DIFaddend}{\DIFOaddend \let\includegraphics\DIFOincludegraphics} 
\DeclareRobustCommand{\DIFdelbegin}{\DIFOdelbegin \let\includegraphics\DIFdelincludegraphics} 
\DeclareRobustCommand{\DIFdelend}{\DIFOaddend \let\includegraphics\DIFOincludegraphics} 
\LetLtxMacro{\DIFOaddbeginFL}{\DIFaddbeginFL} 
\LetLtxMacro{\DIFOaddendFL}{\DIFaddendFL} 
\LetLtxMacro{\DIFOdelbeginFL}{\DIFdelbeginFL} 
\LetLtxMacro{\DIFOdelendFL}{\DIFdelendFL} 
\DeclareRobustCommand{\DIFaddbeginFL}{\DIFOaddbeginFL \let\includegraphics\DIFaddincludegraphics} 
\DeclareRobustCommand{\DIFaddendFL}{\DIFOaddendFL \let\includegraphics\DIFOincludegraphics} 
\DeclareRobustCommand{\DIFdelbeginFL}{\DIFOdelbeginFL \let\includegraphics\DIFdelincludegraphics} 
\DeclareRobustCommand{\DIFdelendFL}{\DIFOaddendFL \let\includegraphics\DIFOincludegraphics} 

\begin{document}
\begin{sloppypar}
\begin{center}
{\Large\bf Covering dense $r$-graphs by perfect matchings}
\\[20pt]
{Jingchao\ Lai, Danyi\ Li$^*$, Weigen\ Yan}\footnote{*Corresponding author.\newline\hspace*{5mm}{Supported by NSFC (Grant Nos. 12571366, 12071180, 12601665).} \newline\hspace*{5mm}{Email addresses: jclai1999@163.com (J. Lai), danyilee1007@163.com (D. Li), weigenyan@263.net (W. Yan).} }
\\[10pt]
\footnotesize { 
School of Science, Jimei University, Xiamen 361021, China}\\
\end{center}
\title{}

\begin{abstract}
For a matching-covered graph $G$, let $\tau(G)$ denote the minimum number of perfect matchings needed to cover all edges of $G$. An $r$-regular graph $G$ is said to be an $r$-graph if $|\partial(X)|\geq r$ for each odd set $X\subseteq V(G)$, where $\partial(X)$ denotes the set of edges with precisely one end in $X$.
The generalized Berge conjecture, which was proposed by Mazzuoccolo (J. Graph Theory, 73 (2013) 377-385) and  equivalent to the generalized Berge-Fulkerson Conjecture posed by Seymour (Proc.
London Math. Soc., 38 (1979) 423-460), states that $\tau(G)\leq 2r-1$ for every $r$-graph $G$, while the 1-factorization conjecture, first explicitly formulated by Chetwynd and Hilton (Proc. London Math. Soc., 50 (1985) 193-206), states that $\tau(G)=r$ for every $r$-graph $G$ of order $2n$ with $r\ge n$. These conjectures are still open.
In this paper, we obtain the following two main results:

1. Let $G$ be an $r$-graph of order $2n$. If $r\ge n$, then $\tau(G)\leq 2r-1$, which implies that the generalized Berge conjecture holds for $r\ge n$, and a linear upper bound of $\tau(G)$ in the 1-factorization conjecture is obtained.

2. Let $G$ be an $r$-graph of order $2n$. If $r\ge 2n/3$, then $\tau(G)\leq 4r$, which implies that a linear upper bound of $\tau(G)$ in the generalized Berge conjecture is obtained if $r\geq 2n/3$.

Our main approach combines equitable proper edge-colorings, matching extension techniques, and applications of Edmonds's perfect matching polytope theorem.
\end{abstract}

\maketitle
\section{Introduction}
All graphs considered in this paper are simple unless otherwise specified. For general graph-theoretic terminology, we follow Bondy and Murty \cite{BM08}.

Let $G$ be a graph. Denote by $V(G)$, $E(G)$, $\delta(G)$ and $\Delta(G)$ the vertex set, the edge set, the minimum degree and the maximum degree of $G$, respectively. For any vertex $v$ in $G$, let $d_G(v)$ be the degree of vertex $v$ and $N_G(v)$ be the set of vertices that are adjacent to $v$. When no confusion arises, we write $d(v)$ and $N(v)$ instead of $d_G(v)$ and $N_G(v)$, respectively. For a vertex subset $X\subseteq V(G)$, set $N(X)=\cup_{x\in X}N(x)$, and the subgraph of $G$ induced by $X$ is denoted by $G[X]$, and $\partial_G(X)$ denotes the set of edges in $G$ with precisely one end in $X$, and $G-H$ (or $G-V(H)$) denotes the subgraph $G[V(G)\setminus V(H)]$ for a subgraph $H$ of $G$.
For two subsets $V_1$ and $V_2$ of $V(G)$, where $V_1\cap V_2=\emptyset$, let $E_G(V_1,V_2)$ denote the set of edges in $G$ with one end in $V_1$ and the other end in $V_2$, and let $e_G(V_1,V_2)=|E_G(V_1,V_2)|$. For convenience, given a subgraph $H$ of $G$, we write $E_H(V_1,V_2) := E(H)\cap E_G(V_1,V_2)$ and $e_H(V_1,V_2):= |E_H(V_1,V_2)|$. A {\em perfect matching} of a graph $G$ is a matching covering all vertices of $G$, and a {\em near-perfect matching} of a graph $G$ is a matching that covers all vertices of $G$ except exactly one vertex.  A graph is {\em matching-covered} if every edge belongs to a perfect matching. An edge in a graph $G$ is called a {\em bridge} if its removal increases the number of components of $G$, and a graph is {\em bridgeless} if it contains no bridge. For a matching-covered graph $G$, let $\tau(G)$ denote the minimum number of perfect matchings needed to cover all edges of $G$.

A {\em proper $k$-edge-coloring} of $G$ is a mapping $f: E(G) \to \{1,\dots,k\}$ such that any two adjacent edges receive distinct colors. That is, $ f(e_1) \neq f(e_2)$ for any adjacent $e_1, e_2 \in E(G)$. The chromatic index of $G$, denoted by $\chi'(G)$, is the smallest integer $k$ such that $G$ admits a proper $k$-edge-coloring.  The set of all edges assigned the same color in a proper edge-coloring is called a {\em color class}. Obviously, by the definition of proper edge-coloring, every color class is a matching of $G$, and a proper $k$-edge-coloring corresponds to a decomposition of $E(G)$ into $k$ pairwise edge-disjoint matchings, some of which may be empty.
A proper edge-coloring of a graph $G$ is said to be {\em equitable} if the sizes of any two color classes differ by at most one.  Vizing's well-known theorem \cite{Vizing} says that $\chi'(G)$ is either $\Delta(G)$ or $\Delta(G)+1$.
For a positive integer $n$, define $[n]=\{1,\ldots,n\}$.

A well-known conjecture, proposed by Berge, states as follows:

\begin{conjecture}[Berge Conjecture]\label{B-conj}
Every bridgeless 3-regular graph $G$ has five perfect matchings such that each edge of $G$ is contained in at least one of them. That is, $\tau(G)\leq 5$ for every bridgeless 3-regular graph $G$.
\end{conjecture}

Another well-known conjecture, which was attributed to Berge in \cite{S1979} but was first published by Fulkerson in \cite{F1971},  states as follows:

\begin{conjecture}[Berge-Fulkerson Conjecture]\label{BF-conj}
Every bridgeless 3-regular graph $G$ has six perfect matchings such that each edge of $G$ is contained in precisely two of them.
\end{conjecture}

In 2011, Mazzuoccolo \cite{Maz2011} proved that the previous two conjectures are equivalent.

An $r$-regular graph is said to be an $r$-graph if $|\partial(X)|\geq r$ for each odd set $X\subseteq V(G)$. It is easy to see that a 3-regular graph is a 3-graph if and only if it is bridgeless. The notion of $r$-graphs was first introduced by Seymour, who also showed that an $r$-graph is always matching-covered \cite{S1979}. Hence, it is natural to ask analogous questions about perfect matching covers on $r$-graphs, in line with the setting of bridgeless 3-regular graphs. In particular, the two conjectures have been generalized to the class of $r$-graphs.

In 1979, Seymour \cite{S1979} proposed the generalized Berge-Fulkerson conjecture:

\begin{conjecture}[Generalized Berge-Fulkerson Conjecture]\label{GBF-conj}
Every $r$-graph $G$ has $2r$ perfect matchings such that each edge of $G$ is contained in exactly two of them.
\end{conjecture}

It is straightforward to see that Conjecture \ref{GBF-conj} implies the following generalization of Conjecture \ref{B-conj}, which was first proposed by Mazzuoccolo in \cite{Maz2013JGT}.

\begin{conjecture}[Generalized Berge Conjecture]\label{GB-conj}
Every $r$-graph $G$ has $2r-1$ perfect matchings such that each edge of $G$ is contained in at least
one of them. That is, $\tau(G)\leq 2r-1$ for every $r$-graph $G$.
\end{conjecture}

In the same paper, Mazzuoccolo \cite{Maz2013JGT} established the equivalence between the generalized Berge-Fulkerson conjecture and the generalized Berge conjecture using an argument analogous to the one he employed for  $3$-graphs. Jin and Steffen \cite{JS2020} established lower bounds on the fraction of edges in an $r$-graph that can be covered by a fixed number of perfect matchings. Following this line of work, Chen and Tokar \cite{CT2023} derived recursive lower bounds for these covering ratios. To date, progress on Conjecture \ref{GB-conj} is very limited, and the problem  is still wide open.

The following conjecture is known as the 1-factorization conjecture.

\begin{conjecture}[\upshape\cite{CH1985}]
Let $G$ be an $r$-graph of order $2n$. If $r\geq n$, then $\tau(G)=r$.
\end{conjecture}
\begin{remark}
It is worth noting that the 1-factorization conjecture is usually formulated for $r$-regular graphs rather than $r$-graphs. However, when $r\geq n$, it is easy to verify that every $r$-regular graph of order $2n$ is an $r$-graph (see Lemma \ref{rr}). Thus, the two formulations are equivalent in this setting.

\end{remark}

The 1-factorization conjecture was first stated explicitly by Chetwynd and Hilton \cite{CH1985}, who also acknowledged that, according to Dirac, the conjecture had been discussed as early as the 1950s. In \cite{CH1985}, Chetwynd and Hilton showed that the conjecture holds for $r\geq12n/7 $. Later, the same authors \cite{CH1989} and, independently, Niessen and Volkmann \cite{NV1990} improved this bound to $(\sqrt{7}-1)n$. In 2016, Csaba et al. \cite{Birminham} confirmed the conjecture for sufficiently large $n$. However, without the assumption that $n$ is sufficiently large, the conjecture remains open in general in the range $n\leq r< (\sqrt{7}-1)n$.

In \cite{CT2023}, Chen and Tokar also obtained the following result.

\begin{theorem}[\upshape\cite{CT2023}]\label{CT2023}
Let $G$ be an $r$-graph of order $2n$, where $r\ge 3$. If $r\geq n$, then $\tau(G)\leq \lceil3r\mathrm{log}r\rceil$.
\end{theorem}

In this paper, we obtain mainly the following two results.

\begin{theorem}\label{main-1}
Let $G$ be an $r$-graph of order $2n$. If $r\ge n$, then $\tau(G)\leq 2r-1$.
Equivalently, the generalized Berge conjecture holds for $r\ge n$.
\end{theorem}

\begin{remark}
The 1-factorization conjecture predicts that $\tau(G)=r$ for $r\ge n$. Although Theorem \ref{main-1} does not attain this optimal bound, it establishes the generalized Berge conjecture throughout the same range.
Moreover, it yields a linear upper bound $\tau(G)\le 2r-1$ for the 1-factorization conjecture, which improves the one in Theorem \ref{CT2023}.
\end{remark}

\begin{theorem}\label{main-2}
Let $G$ be an $r$-graph of order $2n$. If $r\ge 2n/3$, then $\tau(G)\leq 4r$.
\end{theorem}

\begin{remark}
The existence of an absolute constant $C$ satisfying $\tau(G)\le Cr$ for every $r$-graph $G$ remains open. The boundedness problem is unresolved even when $r=3$. Theorem \ref{main-2} establishes such a linear bound with $C=4$ for the class of $r$-graphs of order $2n$ satisfying $r\ge 2n/3$, providing a partial result toward this weaker form of the generalized Berge conjecture.
\end{remark}

We first briefly describe the proof strategy. The proofs of our main results combine equitable proper edge-colorings and matching extension techniques with applications of Edmonds's perfect matching polytope theorem \cite{E1965}. For a graph $G$, we first obtain an equitable proper edge-coloring (that is, a decomposition of $E(G)$ into pairwise edge-disjoint matchings where any two matchings differ in size by at most one), and then extend these matchings to perfect matchings while controlling the number of discarded edges. Specifically, for Theorem \ref{main-1}, we use the degree condition $r\ge n$ to extend matchings with few modifications. To deal with the edges that remain uncovered, we construct a suitable fractional perfect matching and apply Edmonds's theorem \cite{E1965} to find additional perfect matchings that cover multiple uncovered edges simultaneously; for Theorem \ref{main-2}, we establish several structural lemmas based on local matching exchanges, which allow us to extend matchings while keeping the number of discarded edges bounded. We then exploit the uniform fractional perfect matching of an $r$-graph to iteratively cover a guaranteed proportion of the remaining edges, leading to the desired upper bound.

The remainder of this paper is organized as follows. In Section 2, we will briefly recall the
Perfect Matching Polytope Theorem of Edmonds \cite{E1965}, and establish several preliminary results for Theorem \ref{main-1}. In Section 3, we develop several structural lemmas for matchings in dense $r$-graphs, which play a crucial role in the proof of Theorem \ref{main-2}. In Section 4, we present the proofs of Theorems \ref{main-1} and \ref{main-2}.

\section{Preliminaries for Theorem \ref{main-1}}
Let $G$ be a graph. An {\em edge-cut} in $G$ is a set of edges, $C\subseteq E(G)$, such that $C=\partial(X)$ for some $X\subseteq V(G)$. In particular, we call $C=\partial(X)$ {\em odd} if $|X|$ is odd. An odd-cut $C=\partial(X)$ is {\em nontrivial} if $|X|\neq 1$ and $|V(G)\setminus X|\neq 1$.

Let $\mathbf{w}$ be a vector in $\mathbb{R}^{E(G)}$. The entry of $\mathbf{w}$ corresponding to an edge $e$ is denoted by $\mathbf{w}(e)$, and for a subset $A\subseteq E(G)$, we define $\mathbf{w}(A)=\sum\limits_{e\in A}\mathbf{w}(e)$, and the {\em characteristic vector} $\chi^A \in \mathbb{R}^{E(G)}$ of $A$ by $\chi^A(e)=1$ if $e\in A$, and $\chi^A(e)=0$ otherwise.
The vector $\mathbf{w}$ is said to be a {\em fractional perfect matching} if it
satisfies the following:

(a) $0\leq \mathbf{w}(e)\leq 1$ for each edge $e\in E(G)$;

(b) $\mathbf{w}(\partial(v))=1$ for each vertex $v\in V(G)$;

(c) $\mathbf{w}(\partial(X))\geq1$ for each odd set $X\subseteq V(G)$.

The {\em perfect matching polytope} of $G$, denoted by $\mathcal{F}(G)$, is the set of all fractional perfect matchings of $G$. It is easy to see that $\mathbf{w}=(\frac{1}{r},\frac{1}{r},\ldots,\frac{1}{r})\in \mathcal{F}(G)$  for any $r$-graph $G$.

\begin{theorem}[\upshape Perfect Matching Polytope Theorem \cite{E1965}]\label{E1965}
For any graph $G$, the set $\mathcal{F}(G)$ coincides with the convex hull of the characteristic vectors of perfect matchings of $G$.
\end{theorem}

\begin{lemma}\label{Kn}
Let $K_n$ be a complete graph with $V(K_n)=\{u_1,\ldots,u_n\}$, where $n\geq 5$ is odd, and let $e\in E(K_n)$. \\
{\em (1).} There exist $n$ pairwise edge-disjoint near-perfect matchings $F_{u_1},\ldots,F_{u_n}$ such that
$$E(K_n)=\bigcup_{i=1}^{n}F_{u_i},$$
where $u_i$ is the unique uncovered vertex of $F_{u_i}$.\\
{\em (2).} There exist $n$ near-perfect matchings $T_{u_1},\ldots,T_{u_n}$ such that
$$\left|E(K_n-e)\setminus \bigcup_{i=1}^nT_{u_i}\right|=1,$$
where $u_i$ is the unique uncovered vertex of $T_{u_i}$.
\end{lemma}
\begin{proof}
Adding a new vertex $v$ adjacent to every vertex of $K_n$, we obtain $K_{n+1}$. Since $n$ is odd, $K_{n+1}$ can be decomposed into $n$ pairwise edge-disjoint perfect matchings $M_1,\ldots,M_n$, where $vu_i\in M_i$ for each $i\in [n]$. Let $F_{u_i}=M_i\setminus\{vu_i\}.$
Then $F_{u_1},\ldots,F_{u_n}$ are $n$ pairwise edge-disjoint near-perfect matchings of $K_n$ with unique uncovered vertices $u_i$, and
$$E(K_n)=\bigcup_{i=1}^{n}F_{u_i}.$$
This proves the statement (1).

For the statement (2), since $n\ge 5$ is odd, we have $|F_{u_i}|\ge 2$ for every $i\in [n]$. Without loss of generality, we may assume that $e=u_1u_2\in F_{u_n}$. Since $n\ge5$, $F_{u_n}$ contains another edge, say $u_3u_4$. Let
$$T_{u_n}=(F_{u_n}\setminus \{u_1u_2,u_3u_4\})\cup \{u_1u_3, u_2u_4\}$$
and $T_{u_i}=F_{u_i}$ for each $i\in [n-1]$. Then $T_{u_i}$ is a near-perfect matching of $K_n-e$ with unique uncovered vertex $u_i$ for each $i\in [n]$.  Since $\{u_1u_2,u_3u_4\}\subseteq F_{u_n}$, $\{u_1u_3,u_2u_4\}\subseteq\bigcup_{i=1}^{n-1}F_{u_i}=\bigcup_{i=1}^{n-1}T_{u_i}$. Hence
$$\bigcup_{i=1}^n T_{u_i}=E(K_n)\setminus\{u_1u_2,u_3u_4\},$$
and therefore
$$\left|E(K_n-e)\setminus\bigcup_{i=1}^n T_{u_i}\right|=|\{u_3u_4\}|=1.$$
This proves the statement (2). The lemma thus holds.
\end{proof}

\begin{lemma}\label{n,n+2}
Let $G$ be an $n$-graph of order $2n$, where $n\ge 5$. If $G$ has a nontrivial odd-cut of size $n$ or $n+2$, then $\tau(G)\leq 2n-1.$
\end{lemma}
\begin{proof}
Let $V=V(G)$ and let $C$ be a nontrivial odd cut of $G$ such that $c:=|C|\leq n+2$ and $C=\partial(S)$ for some $S\subseteq V(G)$. Without loss of generality, we assume that $s:=|S|\leq |V\setminus S|$. Then $3\leq s\leq n$. Note that
\begin{equation*}\label{f(s)}
c=\sum_{u\in S}d(u)-2|E(G[S])|\ge ns-2\binom{s}{2}=-s^2+(n+1)s:=f(s).
\end{equation*}
Note that $f(s)$ is a concave quadratic with axis of symmetry $s=(n+1)/2$. If $3\le s\le n-1$, then
$$f(s)\ge f(n-1)=2n-2>n+2\ge c,$$
where the second inequality follows from $n\ge 5$, a contradiction. Thus $n=s$ is odd. So $|S|=|V\setminus S|=n$. Let $S=\{u_1,\ldots,u_n\}$ and $V\setminus S=\{v_1,\ldots,v_n\}$.

If $c=n$, then
$$n=c=\sum_{u\in S}d_G(u)-2|E(G[S])|=n^2-2|E(G[S])|.$$
That is, $|E(G[S])|=\binom{n}{2}$ and hence $G[S]\cong K_n$. Similarly, $G[V\setminus S]\cong K_n$.
Since $G$ is $n$-regular, $e_G(\{u_i\},V\setminus S)=1$ and $e_G(\{v_i\},S)=1$ for $i\in [n]$, which implies that $C$ is a perfect matching. Without loss of generality, we assume that $u_iv_i\in E(G)$ for each $i\in [n]$. By (1) of Lemma \ref{Kn}, there exist $2n$ pairwise edge-disjoint near-perfect matchings $F_{u_1},\ldots,F_{u_n},F_{v_1},\ldots, F_{v_n}$ such that
$$E(G[S])=F_{u_1}\cup\cdots \cup F_{u_n}\;\;\text{and}\;\;E(G[V\setminus S])=F_{v_1}\cup\cdots \cup F_{v_n},$$
where $u_i,v_i$ are the unique uncovered vertices of $F_{u_i},F_{v_i}$, respectively. For each $i\in [n]$, let
$$P_i=F_{u_i}\cup \{u_iv_i\}\cup F_{v_i}.$$
Then $P_1,\ldots,P_n$ are $n$ pairwise edge-disjoint perfect matchings such that $E(G)=\bigcup_{i=1}^nP_i$. Hence $\tau(G)=n\leq 2n-1$.

If $c=n+2$, then
$$n+2=c=\sum_{u\in S}d_G(u)-2|E(G[S])|=n^2-2|E(G[S])|.$$
That is, $|E(G[S])|=\binom{n}{2}-1$. Similarly, $|E(G[V\setminus S])|=\binom{n}{2}-1$. So $G[S]\cong G[V\setminus S]\cong K_n-e$, where $e\in E(K_n)$. Moreover, $e_G(\{u_i\},V\setminus S)\ge1$ and $e_G(\{v_i\},S)\ge 1$ for $i\in [n]$. Thus,
\begin{equation}\label{eC}
\bigcup_{e=u_iv_j\in C}\{u_i,v_j\}=\bigcup_{i\in [n]}\{u_i,v_i\}.
\end{equation}
By (2) of Lemma \ref{Kn}, there exist $2n$ near-perfect matchings $T_{u_1},\ldots,T_{u_n},T_{v_1},\ldots,T_{v_n}$ such that
\begin{equation}\label{SVS}
\left|E(G[S])\setminus \bigcup_{i=1}^nT_{u_i}\right|=1\;\;\text{and}\;\;
\left|E(G[V\setminus S])\setminus \bigcup_{i=1}^nT_{v_i}\right|=1,
\end{equation}
where $u_i,v_i$ are the unique uncovered vertices of $T_{u_i},T_{v_i}$, respectively. For every $e=u_iv_j\in C$, where $i,j\in [n]$, let
$$P_e=T_{u_i}\cup \{u_iv_j\}\cup T_{v_j}.$$
Then $P_e$ is a perfect matching of $G$. By \eqref{eC}, we have
\begin{align}\label{Pe}
\bigcup_{e=u_iv_j\in C}P_e&=\bigcup_{e=u_iv_j\in C}\left(T_{u_i}\cup \{u_iv_j\}\cup T_{v_j}\right)=\left(\bigcup_{i=1}^n T_{u_i}\right)\cup\left(\bigcup_{j=1}^n T_{v_j}\right)\cup C.
\end{align}
Note that $E(G)=E(G[S])\cup E(G[V\setminus S])\cup C$. Hence by \eqref{SVS} and \eqref{Pe}, we have
$$\left|E(G)\setminus \bigcup_{e=u_iv_j\in C}P_e\right|=2.$$
Since $G$ is an $n$-graph, $G$ is matching-covered, and hence $E(G)\setminus \bigcup_{e=u_iv_j\in C}P_e$ can be covered by the union of at most 2 perfect matchings. Therefore,
$$\tau(G)\le |C|+2=n+4\le 2n-1,$$
where the last inequality follows from $n\ge 5$. The result thus holds.
\end{proof}

\begin{lemma}\label{|MP|-1}
Let $G$ be a graph of order $2n$ that has a perfect matching. If $d(x)+d(y)\ge 2n-1$ for any two nonadjacent vertices $x$ and $y$ of $G$, then for any matching $M$ of size $n-d$, where $d\ge 0$ is an integer, there exists a perfect matching $P$ such that
$$|M\setminus P|\le d.$$
\end{lemma}
\begin{proof}
We argue by induction on $d$. If $d=0$, then it is trivial. We now assume that $d\ge1$ and the conclusion holds for any matching of size at least $n-d+1$. Let $M=\{u_1v_1,\ldots,u_{n-d}v_{n-d}\}$ and $U=V(G)\backslash V(M)$. Then $|U|=2d\ge 2$.
If there exist $x,y\in U$ satisfying $xy\in E(G)$, then let $M'=M\cup\{xy\}$. Thus $M'$ is a matching of size $n-(d-1)$ and $|M\backslash M'|=0$. By induction on $d-1$, there exists a perfect matching $P$ satisfying $|M'\setminus P|\le d-1$. Then
$$|M\setminus P|\leq |M\setminus M'|+|M'\setminus P|\le 0+d-1<d.$$

It remains to consider the case that $U$ is independent. Choose distinct vertices $x,y\in U$. There exists an index $i\in [n-d]$ such that $e_G(\{x,y\},\{u_i,v_i\})\ge 3$; otherwise, $2n-1\le d(x)+d(y)\le 2(n-d)\le 2n-2$, a contradiction. Without loss of generality, assume that $e_G(\{x,y\},\{u_1,v_1\})\ge 3$, so we have $xu_1,yv_1\in E(G)$ or $xv_1,yu_1\in E(G)$, say $xu_1,yv_1\in E(G)$. Let $M'=(M \setminus \{u_1v_1\})\cup\{xu_1,yv_1\}$. Then $M'$ is a matching of size $n-(d-1)$ and $|M\backslash M'|=1$. By induction on $d-1$, there exists a perfect matching $P$ satisfying $|M'\setminus P|\le d-1$. Then
$$|M\setminus P|\leq |M\setminus M'|+|M'\setminus P|\le 1+d-1=d.$$
The lemma thus holds.
\end{proof}

\begin{theorem}\label{chi}
Let $G$ be a matching-covered graph of order $2n$ and size $m$. If $d(x)+d(y)\ge 2n-1$ for any two nonadjacent vertices $x$ and $y$ of $G$, then $\tau(G)\le (n+1)\chi'(G)-m$.
\end{theorem}
\begin{proof}
Write $\chi'=\chi'(G)$. By the definition of $\chi'$, we have $E(G)=M_1\cup\ldots \cup M_{\chi'}$, where $|M_i|=n-d_i$. Hence
\begin{equation}\label{d_i}
m=|E(G)|=\sum_{i=1}^{\chi'} |M_i|=\chi' n-\sum_{i=1}^{\chi'} d_i.
\end{equation}
By Lemma \ref{|MP|-1}, there exists a perfect matching $P_i$ satisfying $|M_i\setminus P_i|\le d_i$ for each $i\in [\chi']$. Let
$$R=E(G)\setminus \bigcup_{i=1}^{\chi'} P_i.$$
For every edge $e\in R$, there exists a unique index $i\in[\chi']$ such that $e\in M_i$. Since $e\notin P_i$, $e\in M_i\setminus P_i$, which implies that $R\subseteq \bigcup_{i=1}^{\chi'} (M_i\setminus P_i)$. Thus,
\begin{equation}\label{R-1}
|R|\le \left|\bigcup_{i=1}^{\chi'} (M_i\setminus P_i)\right|\le \sum_{i=1}^{\chi'}|M_i\setminus P_i|\le \sum_{i=1}^{\chi'} d_i.
\end{equation}
By \eqref{d_i} and \eqref{R-1}, we have $|R|\le \chi' n-m$. Since $G$ is matching-covered, $R$ can be covered by the union of at most $\chi' n-m$ perfect matchings. Consequently,
$$\tau(G)\le \chi' +(\chi' n-m)=(n+1)\chi'-m.$$
The result thus follows.
\end{proof}


\begin{corollary}
Let $G$ be an $r$-graph of order $2n$. If $r\ge n$, then $\tau(G)\le r+n+1$. In particular, $\tau(G)\le 2r-1$ if $r\ge n+2$. That is, the generalized Berge conjecture holds for $r\ge n+2$.
\end{corollary}
\begin{proof}
Since $G$ is an $r$-graph, $G$ is matching-covered. Moreover, for any two nonadjacent vertices $x$ and $y$ of $G$, we have $d(x)+d(y)=2r> 2n-1$. By Vizing's theorem \cite{Vizing}, $\chi'(G)\le r+1$. Hence, by Theorem \ref{chi},
$$\tau(G)\le (n+1)\chi'(G)-m\le (n+1)(r+1)-m=r+n+1,$$
where the last equality follows from $m=rn$, since $G$ is $r$-regular.
\end{proof}

\section{Matching extension lemmas}
Let $M$ be a matching of a graph $G$. For each $x\in V(M)$, define $m(x)$ to be the unique vertex such that $xm(x)\in M$. Thus the map $m:V(M)\to V(M)$ is a fixed-point-free involution; that is,
$$
m(x)\ne x
\quad\text{and}\quad
m(m(x))=x
$$
for every $x\in V(M)$. In particular, $m$ is a bijection. For a subset $X\subseteq V(M)$, define $m(X)=\{m(x):x\in X\}$. Then $|m(X)|=|X|$ and $m(X)\subseteq N(X)$. Moreover, $X\cap m(X)\ne\emptyset$ if and only if some edge of $M$ has both endpoints in $X$.

\begin{lemma}\label{9-1}
Let $r\ge 1$ be an integer. Let $G$ be a graph of order $2n\le 3r$ with $\delta(G)\ge r$. Let $M$ be a matching with exactly two uncovered vertices $u,v$, and suppose $N(u)\cap N(v)\neq \emptyset$. If for every $xy\in M$ we have $d(x)=d(y)$, then there exists a perfect matching $P$ such that $|M\backslash P|\le 4$.
\end{lemma}
\begin{proof}
We argue by contradiction. If $uv\in E(G)$, then $P:=M\cup \{uv\}$ is a perfect matching such that $|M\setminus P|=0<4$, a contradiction. Hence, we may assume that $uv\notin E(G)$. So $N(u)\cup N(v)\subseteq V(M)$. Let $A=N(u)\cap N(v)$ and $a=|A|$. So $a\ge1$. Let
$$K=V(G)\setminus\left(\{u,v\}\cup m(N(u)\cup N(v))\right).$$
That is, for each $h\in V(M)$, $h\in K$ if and only if $m(h)\notin N(u)\cup N(v)$. Note that $|m(N(u)\cup N(v))|=|N(u)\cup N(v)|=d(u)+d(v)-a$. Thus,
\begin{equation}\label{K}
|K|=2n-2-(d(u)+d(v)-a)\le 3r-(r+r-a)-2=r+a-2.
\end{equation}

We now present the following four successive claims.

\begin{claim}\label{A,m(A)}
$A\cap m(A)=\emptyset$, $\{u,v\}\cup m(A)$ is independent, $m(A)\cap K=\emptyset$ and $A\subseteq K$.
\end{claim}
\begin{proof}
Note that $A\cap m(A)=\emptyset$ if and only if there is no edge in $M$ with both endpoints in $A$. Otherwise, let $e=w_1w_2$ be such an edge in $M$ and $P=(M\setminus \{w_1w_2\})\cup \{uw_1,vw_2\}$. Then $P$ is a perfect matching and $|M\setminus P|=1<4$, a contradiction. So $A\cap m(A)=\emptyset$.

We note that $A\cap m(A)=\emptyset$ implies that for any two distinct vertices $w_1,w_2\in A$, $w_1m(w_1)$ and $w_2m(w_2)$ are two distinct edges in $M$.

If there exist distinct vertices $w_1,w_2\in A$ such that $m(w_1)m(w_2)\in E(G)$, then let $P=(M\setminus \{w_1m(w_1),w_2m(w_2)\})\cup \{uw_1,m(w_1)m(w_2),vw_2\}.$ Thus, $P$ is a perfect matching and $|M\setminus P|=2<4$, a contradiction. So $m(A)$ is independent. Furthermore, if there exists $w\in A$ such that $um(w)\in E(G)$ or $vm(w)\in E(G)$, say $um(w)\in E(G)$, then let
$P=(M\setminus \{wm(w)\})\cup \{um(w),vw\}.$ Then $P$ is a perfect matching and $|M\setminus P|=1<4$, a contradiction. So $\{u,v\}\cup m(A)$ is independent.

Note that $m(A)=m(N(u)\cap N(v))\subseteq m(N(u)\cup N(v))$. So $m(A)\cap K=\emptyset$. Since $\{u,v\}\cup m(A)$ is independent, $m(w)\notin N(u)\cup N(v)$ for every $w\in A$. This implies that $w\in K$. So $A\subseteq K$.
\end{proof}

Since $A\subseteq K$, by \eqref{K}, we have
\begin{equation}\label{KA}
|K\setminus A|\leq r+a-2-a=r-2.
\end{equation}
Now, we let
$$B=A\cup \{x\in K\setminus A\;:\;|N(x)\cap m(A)|\ge 2\}.$$
Then $A\subseteq B \subseteq K$, and $|N(x)\cap m(A)|\le 1$ for every $x\in K\setminus B$.
\begin{claim}\label{B,m(B)}
$B\cap m(B)=\emptyset$ and $m(B)$ is independent.
\end{claim}
\begin{proof}
Note that $B\cap m(B)=\emptyset$ if and only if there is no edge in $M$ with both endpoints in $B$.
Otherwise, let $xy\in M$ be such an edge. By Claim \ref{A,m(A)}, $m(A)\cap K=\emptyset$. So $m(A)\cap B=\emptyset$. This implies that $x,y\notin A$, since $y=m(x)$ and $x=m(y)$.  Thus, $x,y\in B\setminus A$. That is,
$$|N(x)\cap m(A)|\ge 2\;\;\text{and}\;\; |N(y)\cap m(A)|\ge 2.$$
Therefore, we can choose distinct vertices $w_1,w_2\in A$ such that $xm(w_1),ym(w_2)\in E(G)$. Since $x,y\notin A$, $\{x,y\}\cap \{w_1,w_2\}=\emptyset$. Thus, $xy,w_1m(w_1),w_2m(w_2)$ are three distinct edges in $M$. That is, $x,y,w_1,m(w_1),w_2,m(w_2)$ are six distinct vertices in $V(M)$. Let
$$P=(M\setminus \{xy,w_1m(w_1),w_2m(w_2)\})\cup \{uw_1,xm(w_1),ym(w_2),vw_2\}.$$
Then $P$ is a perfect matching and $|M\setminus P|=3<4$, a contradiction. So $B\cap m(B)=\emptyset$.

We note that $B\cap m(B)=\emptyset$ implies that for any two distinct vertices $b_1,b_2\in B$, $b_1m(b_1)$ and $b_2m(b_2)$ are two distinct edges in $M$.

Suppose that there exist distinct vertices $b_1,b_2\in B$ such that $m(b_1)m(b_2)\in E(G)$. By Claim \ref{A,m(A)}, $m(A)$ is independent, so $b_1,b_2$ cannot both lie in $A$. It suffices to consider the following two cases.

(i). If $b_1\in A$ and $b_2\in B\setminus A$ (similarly, $b_2\in A$ and $b_1\in B\setminus A$), then $|N(b_2)\cap m(A)|\ge 2$. So we can choose $w\in A$ such that $w\neq b_1$ and $b_2m(w)\in E(G)$. Since $b_2\in B\setminus A$, $w\neq b_2$. So $w,b_1,b_2$ are three distinct vertices in $B$, and hence $wm(w),b_1m(b_1),b_2m(b_2)$ are three distinct edges in $M$. That is, $w,m(w),b_1,m(b_1),b_2,m(b_2)$ are six distinct vertices in $V(M)$. Let
$$P=(M\setminus \{wm(w),b_1m(b_1),b_2m(b_2)\})\cup\{uw,vb_1,b_2m(w),m(b_1)m(b_2)\}.$$
Then $P$ is a perfect matching and $|M\setminus P|=3<4$, a contradiction.

(ii). If $b_1,b_2\in B\setminus A$, then $|N(b_1)\cap m(A)|\ge 2$ and $|N(b_2)\cap m(A)|\ge 2$. Hence we can choose distinct vertices $w_1,w_2\in A$ such that $b_1m(w_1),b_2m(w_2)\in E(G)$. Clearly, $w_1,w_2,b_1,b_2$ are four distinct vertices in $B$. Thus, $w_1m(w_1),w_2m(w_2),b_1m(b_1),b_2m(b_2)$ are four distinct edges in $M$. That is, $w_1,m(w_1),w_2,m(w_2),b_1,m(b_1),b_2,m(b_2)$ are eight distinct vertices in $V(M)$. Let
\begin{align*}
P&=(M\setminus \{w_1m(w_1),w_2m(w_2),b_1m(b_1),b_2m(b_2)\})\\
&\cup \{uw_1,vw_2,m(b_1)m(b_2),b_1m(w_1),b_2m(w_2)\}.
\end{align*}
Then $P$ is a perfect matching and $|M\setminus P|=4$, a contradiction.

So $m(B)$ is independent.
\end{proof}

\begin{claim}\label{N(B)K}
$N(m(B))\subseteq K$.
\end{claim}
\begin{proof}
Suppose that there exists $h\in N(m(B))\cap (\{u,v\}\cup m(N(u)\cup N(v))$). That is, $h\notin K$.  If $h\in \{u,v\}$, then let $b\in B$ such that $hm(b)\in E(G)$. So $m(b)\in N(u)\cup N(v)$, and hence $b\notin K$, which contradicts $b\in B\subseteq K$.

If $h\in m(N(u)\cup N(v))$, then $m(h)\in N(u)\cup N(v)$. Without loss of generality, assume that $m(h)\in N(u)$. So $um(h)\in E(G)$. Let $b\in B$ such that $hm(b)\in E(G)$. Since $h\notin K$, $b\in B\subseteq K$ and $hm(b)\in E(G)$, $h\neq b$ and $h\neq m(b)$. So $hm(h),bm(b)$ are two distinct edges in $M$. That is, $h,m(h),b,m(b)$ are four distinct vertices in $V(M)$.

(i). If $b\in A$, then let
$$P=(M\setminus \{hm(h),bm(b)\})\cup \{um(h),hm(b),vb\}.$$
Then $P$ is a perfect matching and $|M\setminus P|=2<4$, a contradiction.

(ii). If $b\in B\setminus A$, then $|N(b)\cap m(A)|\ge 2$. Hence we can choose $w\in A$ such that $bm(w)\in E(G)$ and $w\neq m(h)$. Since $h\notin K$, $h\notin A$, and hence $w\neq h$. So $wm(w),hm(h)$ are two distinct edges in $M$. Note that $w\neq b$ and $w,b \in B$. By Claim \ref{B,m(B)}, $B\cap m(B)=\emptyset$, so $wm(w),bm(b)$ are two distinct edges in $M$. Consequently, $wm(w),bm(b),hm(h)$ are three distinct edges in $M$. That is, $w,m(w),b,m(b),h,m(h)$ are six distinct vertices in $V(M)$. Let
$$P=(M\setminus \{wm(w),bm(b),hm(h)\})\cup \{um(h), hm(b),bm(w),vw\}.$$
Then $P$ is a perfect matching and $|M\setminus P|=3<4$, a contradiction.

Therefore, $N(m(B))\cap (\{u,v\}\cup m(N(u)\cup N(v)))=\emptyset$. That is, $N(m(B))\subseteq K$.
\end{proof}

\begin{claim}\label{KBKB}
For every $x\in K\setminus B$, we have
$$|N(m(x))\cap B|\ge |N(x)\cap m(B)|.$$
\end{claim}
\begin{proof}
If $|N(x)\cap m(B)|=0$, then it is trivial. So we assume that $|N(x)\cap m(B)|\geq1$. Let $b\in B$ such that $xm(b)\in E(G)$. By Claim \ref{B,m(B)}, $m(B)$ is independent, so $x\notin m(B)$. Combining this with $x\in K\setminus B$, we have $m(x)\notin B\cup m(B)$ and $m(x)\notin N(u)\cup N(v)$.

We now show that $N(m(x))\subseteq K\cup \{m(b)\}$. Suppose that there exists $h\in \{u,v\}\cup m(N(u)\cup N(v))$ with $h\neq m(b)$ such that $hm(x)\in E(G)$. Since $m(x)\notin N(u)\cup N(v)$, $h\notin \{u,v\}$. So we only need to consider that $h\in m(N(u)\cup N(v))$. Hence $m(h)\in N(u)\cup N(v)$, say $m(h)\in N(u)$. So $um(h)\in E(G)$. Since $h\notin K$, $x\in K\setminus B$ and $b\in B$, $h,x,b$ are three distinct vertices. Together with $h\neq m(b)$, $h\neq m(x)$ and $x\neq m(b)$, we obtain that $hm(h), xm(x), bm(b)$ are three distinct edges in $M$. That is, $h,m(h),x,m(x),b,m(b)$ are six distinct vertices in $V(M)$.

(i). If $b\in A$, then let
$$P=(M\setminus \{hm(h), xm(x), bm(b)\})\cup \{um(h),vb,xm(b),hm(x)\}.$$
Then $P$ is a perfect matching and $|M\setminus P|=3<4$, a contradiction.

(ii). If $b\in B\setminus A$, then $|N(b)\cap m(A)|\ge 2$. Hence we can choose $w\in A$ such that $w\neq m(h)$ and $bm(w)\in E(G)$. Moreover, $b\neq w$ and $b\neq m(w)$. Recall that $m(x)\notin N(u)\cup N(v)$. So $m(x)\notin A$. Since $w\in A\subseteq B\subseteq K$, $x\in K\setminus B$, $m(x)\notin A$ and $h\notin K$, $w\neq x$, $w\neq m(x)$ and $w\neq h$. Together with $w\neq m(h), b\neq w, b\neq m(w)$, we obtain that $wm(w)$ is an edge in $M$ distinct from each of $hm(h), xm(x), bm(b)$. Let
$$P=(M\setminus \{hm(h), xm(x), bm(b),wm(w)\})\cup \{um(h),vw,bm(w),xm(b),hm(x)\}.$$
Then $P$ is a perfect matching and $|M\setminus P|=4$, a contradiction.

So $N(m(x))\subseteq K\cup \{m(b)\}$. This implies that
\begin{equation}
|N(m(x))\setminus B|\leq |(K\cup \{m(b)\})\setminus B|\leq |K\setminus B|+1.
\end{equation}
Therefore,
\begin{align}
|N(m(x))\cap B|&=|N(m(x))|-|N(m(x))\setminus B| \notag\\
&\ge r-|K\setminus B|-1 \notag\\
&=r-(|K\setminus A|-|B\setminus A|)-1 \notag\\
&=r-|K\setminus A|+ |B\setminus A|-1 \notag\\
&\ge|B\setminus A|+1, \label{m(x)B}
\end{align}
where the second equality and the second inequality follow from $A\subseteq B\subseteq K$ and \eqref{KA}, respectively. Recall that $B=A\cup \{x\in K\setminus A:|N(x)\cap m(A)|\ge 2\}$. So $|N(x)\cap m(A)|\le 1$ for every $x\in K\setminus B$. Therefore, for $x\in K\setminus B$, we have
\begin{equation}\label{xm(B)}
|N(x)\cap m(B)|=|N(x)\cap m(A)|+|N(x)\cap m(B\setminus A)|\le 1+|m(B\setminus A)|=1+|B\setminus A|.
\end{equation}
So by \eqref{m(x)B} and \eqref{xm(B)}, we have $|N(m(x))\cap B|\ge |N(x)\cap m(B)|$.
\end{proof}

Note that for every $xy\in M$, we have $d(x)=d(y)$. Hence
\begin{equation}\label{one}
\sum_{b\in B}d(b)=\sum_{m(b)\in m(B)}d(m(b)).
\end{equation}
On the other hand,
\begin{equation}\label{two}
\sum_{b\in B} d(b)= \sum_{b\in B} |N(b)\cap m(B)| + \sum_{b\in B} |N(b)\setminus m(B)|.
\end{equation}
Moreover, by Claim \ref{N(B)K}, $N(m(B))\subseteq K$, so we have
\begin{align}
\sum_{m(b)\in m(B)} d(m(b)) &= \sum_{m(b)\in m(B)} |N(m(b))\cap B| + \sum_{m(b)\in m(B)} |N(m(b))\cap (K\setminus B)| \notag \\
&=\sum_{b\in B} |N(b)\cap m(B)| + \sum_{x\in K\setminus B} |N(x)\cap m(B)|. \label{three}
\end{align}
Thus, by \eqref{one}-\eqref{three}, we have
\begin{equation}\label{four}
\sum_{b\in B} |N(b)\setminus m(B)|=\sum_{x\in K\setminus B} |N(x)\cap m(B)|.
\end{equation}
Recall that $\{u,v\}\cap m(K)=\emptyset$ and $B\subseteq K$. Hence
\begin{equation*}\label{uvm(KB)}
\{u,v\}\cup m(K\setminus B)=(\{u,v\}\cup m(K))\setminus m(B)\subseteq V(G)\setminus m(B).
\end{equation*}
Furthermore,
\begin{align*}
\sum_{b\in B} |N(b)\setminus m(B)|&= \sum_{b\in B} \big|N(b)\cap \big(V(G)\setminus m(B)\big)\big| \\
&\ge \sum_{b\in B} \big|N(b)\cap \big(\{u,v\}\cup m(K\setminus B)\big)\big| \\
&= \sum_{b\in B} \big|N(b)\cap \{u,v\}\big| + \sum_{b\in B} \big|N(b)\cap m(K\setminus B)\big|\\
&\ge\sum_{b\in A} |N(b)\cap \{u,v\}| + \sum_{b\in B} |N(b)\cap m(K\setminus B)|\\
&=2a+\sum_{x\in K\setminus B}|N(m(x))\cap B|\\
&>\sum_{x\in K\setminus B}|N(m(x))\cap B|\\
&\ge \sum_{x\in K\setminus B}|N(x)\cap m(B)|,
\end{align*}
where the second inequality, the third equality and the last inequality follow from $A\subseteq B$, $A=N(u)\cap N(v)$ and Claim \ref{KBKB}, respectively,  which contradicts \eqref{four}.

This completes the proof, and the lemma follows.
\end{proof}

\begin{corollary}\label{9-2}
Let $G$ be an $r$-graph of order $2n\le 3r$, and let $M$ be a matching of size $n-2$ in $G$. Then there exists a matching $M'$ of size at least $n-1$ such that $|M\backslash M'|\le 4$.
\end{corollary}
\begin{proof}
Let $U=V(G)\setminus V(M)$. Then $|U|=4$. If $U$ is not independent, then let $f\in E(G[U])$ and $M'=M\cup\{f\}$. Clearly, $M'$ is a matching of size $n-1$ such that $|M\backslash M'|=0\le 4$. So we assume that $U$ is independent. Let $U=\{u,v,s,t\}$. Since $G$ is $r$-regular, there exist two vertices, say $u,v$, in $U$ such that $N(u)\cap N(v)\neq \emptyset$. Otherwise, $|V(G)|\ge 4r+4>3r$, a contradiction. Let
$$G_1=G\cup \{st\}\;\;\text{and}\;\;M_1=M\cup \{st\},$$
where $V(G_1)=V(G)$ and $E(G_1)=E(G)\cup \{st\}$. Then $|V(G_1)|=2n\le 3r$, $\delta(G_1)=r$ and $M_1$ is a matching of size $n-1$ in $G_1$, for which the two uncovered vertices $u,v$ satisfy $N_{G_1}(u)\cap N_{G_1}(v)\neq \emptyset$. Moreover, $d_{G_1}(s)=d_{G_1}(t)=r+1$ for $st\in M_1$, and $d_{G_1}(x)=d_{G_1}(y)=r$ for $xy\in M$. By Lemma \ref{9-1}, there exists a perfect matching $P$ in $G_1$ such that $|M_1\setminus P|\le 4$.

If $st\in P$, then $|M_1\setminus P|=|M\setminus (P\setminus \{st\})|$. Let $M'=P\setminus \{st\}$. Then $M'$ is a matching of size $n-1$ in $G$, and $|M\setminus M'|\le 4$.

If $st\notin P$, then $|M_1\setminus P|=|M\setminus P|+1$. Let $M'=P$. Then $M'$ is a perfect matching in $G$, and $|M\setminus M'|=|M_1\setminus P|-1\le 3$.

The result thus holds.
\end{proof}

\begin{lemma}\label{9-3}
Let $G$ be an $r$-graph of order $2n\le 3r$, and let $M$ be a matching of size $n-1$ in $G$. Then at least one of the following holds.\\
{\em (1).} There exists a perfect matching $P$ such that $|M\backslash P|\le 4$; \\
{\em (2).} There exists a matching $M'$ of size $n-1$ with $|M\setminus M'|\le 2$, such that the two uncovered vertices of $M'$ have a common neighbor.
\end{lemma}
\begin{proof}
We argue by contradiction. Let $V(G)\setminus V(M)=\{u,v\}$. If $uv\in E(G)$, then let $P=M\cup \{uv\}$. Obviously, $P$ is a perfect matching such that $|M\backslash P|=0<4$, a contradiction. So $uv\notin E(G)$. If $N(u)\cap N(v)\neq \emptyset$, then let $M'=M$. Clearly, $M'$ satisfies condition (2), a contradiction. So $N(u)\cap N(v)=\emptyset$. Let
$$X=m(N(u))\subseteq V(M)\;\;\;\text{and}\;\;\;Y=m(N(v))\subseteq V(M).$$
Then $|X|=|Y|=r$ and $u,v\notin X\cup Y$ and
\begin{equation}\label{N(u)N(X)}
N(u)=m(X)\subseteq N(X)\;\;\;\text{and}\;\;\;N(v)=m(Y)\subseteq N(Y).
\end{equation}

We now present the following three successive claims.
\begin{claim}\label{X,Y-1}
$X\cap Y=\emptyset$, $N(u)\cap Y=\emptyset$, $N(v)\cap X=\emptyset$, and $E_G(X,Y)=\emptyset$.
\end{claim}
\begin{proof}
Since $m$ is a bijection, $X\cap Y=m(N(u))\cap m(N(v))=m(N(u)\cap N(v))=\emptyset$.

If there exists $y\in N(u)\cap Y$, then $uy\in E(G)$ and $m(y)\in m(Y)=N(v)$. Hence $vm(y)\in E(G)$. Let
$$P=(M\setminus\{ym(y)\})\cup \{uy,vm(y)\}.$$
Then $P$ is a perfect matching and $|M\setminus P|=1<4$, a contradiction. Thus $N(u)\cap Y=\emptyset$. Similarly, $N(v)\cap X=\emptyset$.

If there exists $xy\in E(G)$ such that $x\in X$ and $y\in Y$, then $m(x)\in m(X)=N(u)$ and $m(y)\in m(Y)=N(v)$. So $m(x)\neq m(y)$ and $um(x),vm(y)\in E(G)$. Since $N(v)\cap X=\emptyset$, we have $x\neq m(y)$. Thus, $xm(x)$ and $ym(y)$ are two distinct edges in $M$. That is, $x,m(x),y,m(y)$ are four distinct vertices in $V(M)$. Let
$$P=(M\setminus \{xm(x),ym(y)\})\cup \{um(x),xy,vm(y)\}.$$
Then $P$ is a perfect matching and $|M\setminus P|=2<4$, a contradiction. So $E_G(X,Y)=\emptyset$.
\end{proof}

\begin{claim}\label{N(X)N(Y)}
$N(X)\cap N(Y)=\emptyset$ and $V(G)=N(X)\cup N(Y)$.
\end{claim}
\begin{proof}
If there exists $z\in N(X)\cap N(Y)$, then $xz,yz\in E(G)$ for some $x\in X$ and some $y\in Y$, and hence $m(x)\in m(X)=N(u)$ and $m(y)\in m(Y)=N(v)$. That is, $um(x),vm(y)\in E(G)$. By Claim \ref{X,Y-1}, $X\cap Y=\emptyset$ and $N(v)\cap X= \emptyset$, so $x\neq y$ and $x\neq m(y)$.  Thus, $xm(x)$ and $ym(y)$ are two distinct edges in $M$. Let
$$M'=(M\setminus \{xm(x),ym(y)\})\cup \{um(x),vm(y)\}.$$
Then $M'$ is a matching of size $n-1$ with $|M\setminus M'|=2$, such that the two uncovered vertices $x,y$ of $M'$ have a common neighbor $z$, a contradiction. So $N(X)\cap N(Y)=\emptyset$.

We next show that $V(G)=N(X)\cup N(Y)$.  If $u\notin N(X)\cup N(Y)$, then $N(u)\cap (X\cup Y)=\emptyset$. Recall that $uv\notin E(G)$, $X\cap Y=\emptyset$ and $v\notin X\cup Y$. Hence $N(u)\subseteq V(G)\setminus (X\cup Y\cup \{v\})$. That is,
$$r=|N(u)|\le 2n-|X|-|Y|-1=2n-2r-1\le 3r-2r-1=r-1,$$
a contradiction. So $u\in N(X)\cup N(Y)$. Similarly, $v\in N(X)\cup N(Y)$.

If there exists $w\in V(G)\setminus \{u,v\}$ such that $w\notin N(X)\cup N(Y)$, then $N(w)\cap (X\cup Y)=\emptyset$. Moreover, $wu\notin E(G)$. Otherwise, $w\in N(u)=m(X)\subseteq N(X)$. Hence $N(w)\subseteq V(G)\setminus (X\cup Y\cup \{u\})$. That is,
$$r=|N(w)|\le 2n-|X|-|Y|-1=2n-2r-1\le 3r-2r-1=r-1,$$
a contradiction. Thus, $V(G)\subseteq N(X)\cup N(Y)$. That is, $V(G)=N(X)\cup N(Y)$.
\end{proof}
By Claims \ref{X,Y-1} and \ref{N(X)N(Y)}, $E_G(X,Y)=\emptyset$ and $V(G)=N(X)\cup N(Y)$. So $X\cap N(Y)=\emptyset$, and hence $X\subseteq N(X)$. Similarly, $Y\subseteq N(Y)$. Moreover, Claim \ref{X,Y-1} gives $N(u)\cap Y=\emptyset$ and $N(v)\cap X=\emptyset$. Consequently, $u\notin N(Y)$ and $v\notin N(X)$.
Combining this with $V(G)=N(X)\cup N(Y)$, we obtain $u\in N(X)$ and $v\in N(Y)$.

\begin{claim}
$E_M(N(X),N(Y))=\emptyset$.
\end{claim}
\begin{proof}
If there exists $ab\in M$ such that $a\in N(X)$ and $b\in N(Y)$, then $ax,by\in E(G)$ for some $x\in X\subseteq N(X)$ and $y\in Y\subseteq N(Y)$. Hence $x\neq y$, $x\neq b$, $y\neq a$ and $m(x)\in m(X)=N(u)$, $m(y)\in m(Y)=N(v)$. By Claim \ref{X,Y-1}, $N(v)\cap X=\emptyset$, so we have $x\neq m(y)$. Thus, $ab,xm(x),ym(y)$ are three distinct edges in $M$. That is, $a,b,x,m(x),y,m(y)$ are six distinct vertices in $V(M)$. Let
$$P=(M\setminus\{ab,xm(x),ym(y)\})\cup \{um(x),ax,by,vm(y)\}.$$
Then $P$ is a perfect matching and $|M\setminus P|=3<4$, a contradiction. So we have $E_M(N(X),N(Y))=\emptyset$.
\end{proof}
Recall that $u\in N(X)$ and $v\in N(Y)$. Since $E_M(N(X),N(Y))=\emptyset$, $M$ pairs every vertex of $N(X)$ internally except $u$. Thus, $|N(X)|$ is odd. Since $G$ is an $r$-graph, we have $|\partial(N(X))|\ge r\geq1$. Consequently, there exists an edge $ab\in E(G)$ with $a\in N(X)$ and $b\in N(Y)$. Obviously, $ab\notin M$. We can see that $a\neq u$, since $b\in N(a)$ and $N(u)=m(X)\subseteq N(X)$. Moreover, $a\notin X$. Otherwise, $b\in N(a)\subseteq N(X)$. Similarly, $b\neq v$ and $b\notin Y$.

Since $E_M(N(X),N(Y))=\emptyset$, $m(a)\in N(X)$ and $m(b)\in N(Y)$. So there exist some $x\in X$ and $y\in Y$ such that $xm(a),ym(b)\in E(G)$. Meanwhile, $m(x)\in m(X)=N(u)$ and $m(y)\in m(Y)=N(v)$. So $um(x),vm(y)\in E(G)$. Since $a\notin X$ and $xm(a)\in E(G)$, $x\neq a$ and $x\neq m(a)$. Thus, $xm(x),am(a)$ are two distinct edges in $M$. That is, $x,m(x),a,m(a)$ are four distinct vertices in $V(M)$. Similarly, $y,m(y),b,m(b)$ are also four distinct vertices in $V(M)$. Note that $x,m(x),a,m(a)\in N(X)$ and $y,m(y),b,m(b)\in N(Y)$. Consequently,
$x,m(x),a,m(a),y,m(y),b,m(b)$ are eight distinct vertices in $V(M)$. Let
$$P=(M\setminus \{xm(x),am(a),ym(y),bm(b)\})\cup \{um(x),xm(a),ab,ym(b),vm(y)\}.$$
Then $P$ is a perfect matching and $|M\setminus P|=4$, a contradiction.

This completes the proof, and the lemma follows.
\end{proof}

\begin{lemma}\label{|MP|-2}
Let $G$ be an $r$-graph of order $2n\le 3r$, and let $M$ be a matching of $G$.\\
{\em (i).} If $|M|=n-1$, then there exists a perfect matching $P$ such that $|M\setminus P|\leq 6$;\\
{\em (ii).} If $|M|=n-2$, then there exists a perfect matching $P$ such that $|M\setminus P|\leq 10$.
\end{lemma}
\begin{proof}
Note that $G$ satisfies the conditions of Lemma \ref{9-3}. For (i), if Lemma \ref{9-3} (1) holds, then (i) follows immediately. Otherwise, we may assume that there exists a matching $M'$ of size $n-1$ with $|M\setminus M'|\le 2$, such that the two uncovered vertices of $M'$ have a common neighbor. Note that $G$ and $M'$ satisfy the conditions of Lemma \ref{9-1}. Hence by Lemma \ref{9-1}, there exists a perfect matching $P$ such that $|M'\setminus P|\le 4$. Thus
$$|M\setminus P|\le |M\setminus M'|+|M'\setminus P|\le 2+4=6.$$

For (ii), by Corollary \ref{9-2}, there exists a matching $M'$ of size at least $n-1$ such that $|M\backslash M'|\le 4$. If $|M'|=n$, i.e., $M'$ is a perfect matching, then it suffices to take $P=M'$. If $|M'|=n-1$, then by (i), there exists a perfect matching $P$ such that $|M'\setminus P|\leq 6$. Thus
$$|M\setminus P|\le |M\setminus M'|+|M'\setminus P|\le 4+6=10.$$
The lemma thus holds.
\end{proof}

\section{Proofs of the main results}
In this section, we complete the proofs of Theorems \ref{main-1} and \ref{main-2}.
We begin with the following lemmas.
\begin{lemma}[\upshape\cite{FF1969}]\label{FF1969}
Let $k$ be a positive integer. If $G$ has a proper $k$-edge-coloring, then $G$ has an equitable proper $k$-edge-coloring.
\end{lemma}

\begin{lemma}\label{rr}
Let $G$ be an $r$-regular graph of order $2n$. If $r\ge n$, then $G$ is an $r$-graph.
\end{lemma}
\begin{proof}
Let $S\subseteq V(G)$ be an odd set. Since $|V(G)|$ is even, $|V(G)\setminus S|$ is also odd, so we may assume that $s:=|S|\le n\le r$. Hence we have
\begin{equation*}
|\partial(S)|= rs-2|E(G[S])|\ge rs- 2\binom{s}{2}=r+(s-1)(r-s) \ge r.
\end{equation*}
Thus $G$ is an $r$-graph.
\end{proof}

Instead of directly proving Theorem \ref{main-1}, we first establish the following result.
\begin{theorem}\label{main-0}
Let $G$ be an $n$-graph of order $2n$. Then $\tau(G)\leq 2n-1$.
\end{theorem}
\begin{proof} In \cite{CH1985}, Chetwynd and Hilton showed that for every connected $d$-regular graph $G$ of order $2n$, $\tau(G)=d$ if $d\ge 2n-5$. For $1\le n\le 4$, we have $2n-5<n$. Hence $\tau(G)=n\le 2n-1$. Therefore, we may assume that $n\ge 5$. For every nontrivial odd-cut $C=\partial(S)$ for some odd set $S\subseteq V(G)$, we have
$$c:=|C|=\sum_{u\in S}d(u)-2|E(G[S])|=n|S|-2|E(G[S])|.$$
Since $|S|$ is odd,  $c\equiv n$ (mod 2).  Since $G$ is an $n$-graph, we have $c\ge n$. If $G$ has a nontrivial odd-cut of size $n$ or $n+2$, then by Lemma \ref{n,n+2}, $\tau(G)\le 2n-1$. Therefore, we may further assume that every nontrivial odd-cut has size at least $n+4$.

By Vizing's theorem \cite{Vizing}, $G$ has a proper $(n+1)$-edge-coloring. Furthermore, by Lemma \ref{FF1969}, $G$ has an equitable proper $(n+1)$-edge-coloring. That is,
$$E(G)=M_1\cup M_2\cup\ldots\cup M_{n+1},$$
where $M_1,\dots,M_{n+1}$ are pairwise edge-disjoint matchings, and the sizes of all matchings differ by at most one. Since $|E(G)| = n^2$,
\[
\frac{|E(G)|}{n+1} = \frac{n^2}{n+1} = n - \frac{n}{n+1}.
\]
Thus, every matching has size $n$ or $n-1$. Let $a$ be the number of matchings of size $n$. Counting edges gives
\[
n^2= an+(n+1-a)(n-1),
\]
hence $a=1$. Therefore, exactly one of the $n+1$ matchings is a perfect matching. Without loss of generality, we assume that $M_1,\ldots,M_n$ are the matchings of size $n-1$ and $M_{n+1}$ is a perfect matching. Let $U_i=V(G)\setminus V(M_i)$ for each $i\in [n]$. Then $|U_i|=2$.

\begin{claim}\label{UUU}
$U_i\cap U_j=\emptyset$ for every pair $i\neq j$, and hence $V(G)=\bigcup_{i=1}^nU_i$.
\end{claim}
\begin{proof}
If there exists $u\in V(G)$ such that $u\in U_i\cap U_j$, then neither $M_i$ nor $M_j$ covers $u$.
Note that $M_1,\ldots,M_{n+1}$ form an edge decomposition of $E(G)$ and $d(u)=n$. Therefore, exactly $n$ of these matchings cover $u$, meaning that exactly one matching fails to cover $u$. This yields a contradiction. So $U_i\cap U_j=\emptyset$ for every pair $i\neq j$. Hence, $|\bigcup_{i=1}^nU_i|=2n$, which implies that $V(G)=\bigcup_{i=1}^nU_i$.
\end{proof}

Clearly, $G$ satisfies the conditions of Lemma \ref{|MP|-1}. Hence, there exists a perfect matching $P_i$ such that $|M_i\setminus P_i|\le 1$ for every $i\in [n]$. Let $P_{n+1}=M_{n+1}$. We define the function $\mathbf{w}: E(G)\rightarrow R$ as follows:
$$\mathbf{w}(e)=\frac{1}{n-2}\left(1-\frac{2}{n+1}\sum_{i=1}^{n+1}\chi^{P_i}(e)\right).$$

\begin{claim}
$\mathbf{w}\in \mathcal{F}(G).$
\end{claim}
\begin{proof}
(a). If $e=uv\in P_i\setminus M_i$, where $i\in [n]$, then $\{u,v\}\cap U_i\neq \emptyset$. Otherwise, $\{u,v\}\subseteq V(M_i)$. Let $x,y$ be two vertices in $V(M_i)$ such that $ux,vy\in M_i$. Obviously, $ux,vy\notin P_i$, since $uv\in P_i$. Thus, $|M_i\setminus P_i|\ge 2$, a contradiction.

We show that $\sum_{i=1}^{n+1}\chi^{P_i}(e)\le 3$ for every $e\in E(G)$. Suppose that there exists $e=uv$ such that $\sum_{i=1}^{n+1}\chi^{P_i}(e)\ge 4$. Note that there is a unique index $j\in [n+1]$ such that $e\in M_j$. Hence, there are at least three indices $i\in [n]$ such that $e\in P_i\setminus M_i$, say $1,2,3$. Hence $|\{u,v\}\cap U_i|\ge 1 $ for $i=1,2,3$. This implies that two sets among $U_1,U_2,U_3$ intersect, yielding a contradiction to Claim \ref{UUU}. So $\sum_{i=1}^{n+1}\chi^{P_i}(e)\le 3$ for every $e\in E(G)$.

Therefore,
$$1-\frac{2}{n+1}\sum_{i=1}^{n+1}\chi^{P_i}(e)\ge1-\frac{6}{n+1}=\frac{n-5}{n+1}\ge 0.$$
So $\mathbf{w}(e)\ge 0$ for every $e\in E(G)$. Clearly, $\mathbf{w}(e)\le 1$ for every $e\in E(G)$.

(b). For every $u\in V(G)$, it is obvious that $|P_i\cap \partial(u)|=1$ for every $i\in[n+1]$. Hence
\begin{align*}
\mathbf{w}(\partial(u))
&= \frac{1}{n-2}\sum_{e\in\partial(u)} \left(1-\frac{2}{n+1}\sum_{i=1}^{n+1}\chi^{P_i}(e)\right) \\
&= \frac{1}{n-2}\left(n-\frac{2}{n+1}\sum_{i=1}^{n+1}\sum_{e\in\partial(u)}\chi^{P_i}(e)\right) \\
&= \frac{1}{n-2}\left(n-\frac{2}{n+1}\sum_{i=1}^{n+1}\big|P_i\cap\partial(u)\big|\right) \\
&= \frac{1}{n-2}\left(n-\frac{2(n+1)}{n+1}\right) = 1.
\end{align*}

(c). For each odd set $S\subseteq V(G)$, let $C=\partial(S)$ and $c=|C|$. If $|S|=1$ or $|V(G)\setminus S|=1$, then by the previous paragraph, we have $\mathbf{w}(C)=1$. It remains to consider the case $|S|\ge 3$ and $|V(G)\setminus S|\ge 3$. Hence $C$ is a nontrivial odd-cut. By assumption, $c\ge n+4$.
Let
$$Q=\{i\in [n]:U_i\cap S\neq \emptyset,U_i\cap (V(G)\setminus S)\neq \emptyset\}$$
and $q=|Q|$. By Claim \ref{UUU}, $U_i\cap U_j=\emptyset$ for every pair $i\neq j$, and hence $V(G)=\bigcup_{i=1}^nU_i$. Moreover, $|U_i|=2$ for each $i\in [n]$. Hence one can see that $q\equiv |S|$ (mod 2). Since $|S|$ is odd, $q\ge 1$. For $i\in[n+1]$, let
$$\triangle_i=|P_i\cap C|-|M_i\cap C|.$$
Clearly, $\triangle_{n+1}=0$. For $i\in [n]$, we have $|P_i\cap M_i|=|M_i|-|M_i\setminus P_i|\ge |M_i|-1= n-2$. So $|P_i\setminus M_i|=|P_i|-|P_i\cap M_i|\le n-(n-2)=2$. Therefore,
\begin{align*}
\triangle_i&=|P_i\cap C|-|M_i\cap C| \\
&=|(P_i\setminus M_i)\cap C|+|(P_i\cap M_i)\cap C|-|M_i\cap C|\\
&\le|(P_i\setminus M_i)\cap C|\\
&\le|P_i\setminus M_i|\le 2.
\end{align*}
For every $i\in Q$, since $U_i$ is separated by $C$, $M_i$ leaves
exactly one vertex of $S$ uncovered. Moreover, $|S|-1$ is even. Thus, $|M_i\cap C|$ is even. Since $|S|$ is odd and $P_i$ is a perfect matching, $|P_i\cap C|$ is odd. So $\triangle_i$ is odd for every $i\in Q$, which implies that $\triangle_i\le 1$ for every $i\in Q$. Hence
\begin{equation}\label{triangle}
\sum_{i=1}^{n+1}\triangle_i=\triangle_{n+1}+\sum_{i\in Q}\triangle_i+\sum_{i\in [n]\setminus Q}\triangle_i
\leq 0+q+2(n-q)=2n-q\le 2n-1.
\end{equation}
Since $M_1,\dots,M_{n+1}$ form an edge-disjoint decomposition of $E(G)$, we have $\sum_{i=1}^{n+1}|M_i\cap C|=|C|=c$. Together with \eqref{triangle}, this yields
\begin{equation}
\sum_{i=1}^{n+1}|P_i\cap C|
= \sum_{i=1}^{n+1}\big(\triangle_i + |M_i\cap C|\big)
= \sum_{i=1}^{n+1}\triangle_i + \sum_{i=1}^{n+1}|M_i\cap C|
\le 2n-1 + c.
\end{equation}
Therefore,
\begin{align*}
\mathbf{w}(C)&= \frac{1}{n-2}\sum_{e\in C}\left(1-\frac{2}{n+1}\sum_{i=1}^{n+1}\chi^{P_i}(e)\right)\\
&=\frac{1}{n-2}\left(c-\frac{2}{n+1}\sum_{i=1}^{n+1}|P_i\cap C|\right)\\
&\ge \frac{1}{n-2}\left[c - \frac{2}{n+1}(2n-1 + c)\right]\\
&= \frac{1}{n-2}\cdot\frac{c(n-1)-4n +2}{n+1}\\
&\ge \frac{(n+4)(n-1)-4n+2}{(n-2)(n+1)}=1.
\end{align*}
By (a), (b) and (c), we have $\mathbf{w}\in \mathcal{F}(G)$. The claim thus follows.
\end{proof}
By Theorem \ref{E1965}, there exist perfect matchings $Q_1,\ldots,Q_m$ such that
\begin{equation}\label{Edmond-1}
\mathbf{w}=\sum_{j=1}^m\lambda_j\chi^{Q_j},
\end{equation}
where $\sum_{j=1}^m\lambda_j=1$ and $\lambda_j>0$. Let
$$R=E(G)\setminus \bigcup_{i=1}^{n+1}P_i.$$
So, for every $e\in R$, $e\notin P_{n+1}=M_{n+1}$, $\sum_{i=1}^{n+1}\chi^{P_i}(e)=0$ and hence $\mathbf{w}(e)=1/(n-2)$; moreover, there exists a unique index $i\in[n]$ such that $e\in M_i$. Since $e\notin P_i$, $e\in M_i\setminus P_i$, which implies that $R\subseteq \bigcup_{i=1}^n (M_i\setminus P_i)$. Thus,
\begin{equation*}
|R|\le \left|\bigcup_{i=1}^n (M_i\setminus P_i)\right|\le \sum_{i=1}^n|M_i\setminus P_i|\le n.
\end{equation*}

If $|R|\le n-2$, then $R$ can be covered by the union of at most $n-2$ perfect matchings, since $G$ is matching-covered. Moreover, $P_1,\ldots,P_{n+1}$ are all perfect matchings. Thus,
$$\tau(G)\leq (n+1)+(n-2)=2n-1.$$

If $n-1\leq |R|\le n$, then $\mathbf{w}(R)=\sum_{e\in R}\mathbf{w}(e)=|R|/(n-2)>1$. By \eqref{Edmond-1},
$$\mathbf{w}(R)=\sum_{e\in R}\mathbf{w}(e)=\sum_{j=1}^m\lambda_j\sum_{e\in R}\chi^{Q_j}(e)=\sum_{j=1}^m\lambda_j|Q_j\cap R|>1,$$
which implies that there exists an index $i\in [m]$, say $i=1$, such that $|Q_1\cap R|\ge 2$, since $\sum_{j=1}^m\lambda_j=1$.  Let $$R_1=R\setminus(Q_1\cap R).$$
Then $R_1\cap Q_1=\emptyset$ and $|R_1|=|R|-|Q_1\cap R|\leq n-2$.

(i). If $|R_1|\leq n-3$, then $R_1$ can be covered by the union of at most $n-3$ perfect matchings, since $G$ is matching-covered. Note that $P_1,\ldots,P_{n+1}$ are all perfect matchings and $Q_1\cap R$ is covered by $Q_1$.  Thus,
$$\tau(G)\leq (n+1)+1+(n-3)=2n-1.$$

(ii). If $|R_1|=n-2$, then $\mathbf{w}(R_1)=\sum_{e\in R_1}\mathbf{w}(e)=|R_1|/(n-2)=1$. By \eqref{Edmond-1},
$$1=\sum_{e\in R_1}\mathbf{w}(e)=\sum_{j=1}^m\lambda_j\sum_{e\in R_1}\chi^{Q_j}(e)=\sum_{j=1}^m\lambda_j|Q_j\cap R_1|=\sum_{j=2}^m\lambda_j|Q_j\cap R_1|,$$
where the last equality follows from $Q_1\cap R_1=\emptyset$. Since $\sum_{j=2}^m\lambda_j=1-\lambda_1<1$, there exists an index $i\in \{2,\ldots,m\}$, say $i=2$, such that $|Q_2\cap R_1|\ge 2$. Choose two distinct edges $e_1,e_2\in Q_2\cap R_1$, and let $R_2=R_1\setminus \{e_1,e_2\}$. So $|R_2|=n-4$. Since $G$ is matching-covered,
$R_2$ can be covered by the union of at most $n-4$ perfect matchings. Thus
$$\tau(G)\le (n+1)+|\{Q_1,Q_2\}|+(n-4)=2n-1.$$
This completes the proof of Theorem \ref{main-0}.
\end{proof}

\vspace{2mm}
\noindent{\bf{Proof of Theorem \ref{main-1}.}} Let $G$ be an $r$-graph of order $2n$, where $r\ge n$.
If $n=1$, then $G=K_2$, so the conclusion is immediate. Hence we may assume that $n\ge2$.

We first successively remove $r-n$ pairwise edge-disjoint perfect matchings from $G$. Indeed, after $j$ perfect matchings have been removed, where $0\le j\le r-n$, the remaining graph is $(r-j)$-regular of order $2n$. Its minimum degree satisfies $r-j\ge n$. By Dirac's theorem \cite{Dirac1952}, the remaining graph has a Hamiltonian cycle and hence a perfect matching. Thus this process can be continued until $r-n$
perfect matchings have been removed. Let $G_0$ be the resulting spanning subgraph. Then $G_0$ is $n$-regular of order $2n$. By Lemma \ref{rr}, $G_0$ is an $n$-graph. By Theorem \ref{main-0}, $\tau(G_0)\le 2n-1$. Hence we have
$$\tau(G)\le \tau(G_0)+r-n\le r+n-1\le 2r-1.$$
This completes the proof of Theorem \ref{main-1}..

\vspace{2mm}
\noindent{\bf{Proof of Theorem \ref{main-2}.}} Let $G$ be an $r$-graph of order $2n\le 3r$.  If $n\le r$, then by Theorem \ref{main-1}, $\tau(G)\le 2r-1<4r$. Thus we only need to consider $r+1\le n\le 3r/2$.

By Vizing's theorem \cite{Vizing}, $G$ has a proper $(r+1)$-edge-coloring. Furthermore, by Lemma \ref{FF1969}, $G$ has an equitable proper $(r+1)$-edge-coloring. That is,
$$E(G)=M_1\cup M_2\cup\ldots\cup M_{r+1},$$
where $M_1,\dots,M_{r+1}$ are pairwise edge-disjoint matchings, and the sizes of all matchings differ by at most one. Since $|E(G)|=rn$,
\[
\frac{|E(G)|}{r+1} = \frac{rn}{r+1} =n-\frac{n}{r+1},
\]
where $1\le n/(r+1)<2$, since $r+1\le n\le 3r/2$. Thus every matching has size $n-1$ or $n-2$. Let $a$ be the number of matchings of size $n-1$. Counting edges gives
\[
rn= a(n-1)+(r+1-a)(n-2),
\]
hence $a=2r+2-n$. Without loss of generality, assume that $M_1,\ldots,M_{2r+2-n}$ are the matchings of size $n-1$, and that $M_{2r+3-n},\ldots,M_{r+1}$ are the matchings of size $n-2$. By Lemma \ref{|MP|-2}, for every $i=1,\ldots,2r+2-n$ and $j=2r+3-n,\ldots,r+1$, there exist perfect matchings $P_i,P_j$ such that
\begin{equation*}
|M_i\setminus P_i|\le 6\;\;\text{and}\;\;|M_j\setminus P_j|\le 10.
\end{equation*}
Let
$$R=E(G)\setminus \bigcup_{i=1}^{r+1}P_i.$$
Then for every $e\in R$, there exists a unique index $i\in[r+1]$ such that $e\in M_i$. Since $e\notin P_i$, $e\in M_i\setminus P_i$, which implies that $R\subseteq \bigcup_{i=1}^{r+1} (M_i\setminus P_i)$. Thus,
\begin{equation}\label{8r+2}
|R|\le \sum_{i=1}^{r+1}|M_i\setminus P_i|\leq 6(2r+2-n)+10(n-r-1)=4n+2r+2\le 8r+2.
\end{equation}

We give the following claim.
\begin{claim}\label{F/r}
For every $F\subseteq E(G)$, there exists a perfect matching $Q$ such that
$$|Q\cap F|\ge \left \lceil\frac{|F|}{r}\right\rceil.$$
\end{claim}
\begin{proof}
Note that $\mathbf{w}=(\frac{1}{r},\frac{1}{r},\ldots,\frac{1}{r})\in \mathcal{F}(G)$. Hence by Theorem \ref{E1965}, there exist perfect matchings $Q_1,\ldots,Q_m$ such that
\begin{equation}\label{Edmond-2}
\mathbf{w}=\sum_{j=1}^m\lambda_j\chi^{Q_j},
\end{equation}
where $\sum_{j=1}^m\lambda_j=1$ and $\lambda_j>0$. Therefore,
$$\mathbf{w}(F)=\sum_{e\in F}\mathbf{w}(e)=\sum_{j=1}^m\lambda_j\sum_{e\in F}\chi^{Q_j}(e)=\sum_{j=1}^m\lambda_j|Q_j\cap F|.$$
On the other hand, $\mathbf{w}(F)=\sum_{e\in F}\mathbf{w}(e)=|F|/r$. That is, $\sum_{j=1}^m\lambda_j|Q_j\cap F|=|F|/r$. Since $\sum_{j=1}^m\lambda_j=1$, there exists an index $j\in [m]$ such that $|Q_j\cap F|\ge |F|/r$. That is, $|Q_j\cap F|\ge \lceil|F|/r\rceil$, since $|Q_j\cap F|$ is an integer.
\end{proof}
By Claim \ref{F/r}, there exists a perfect matching $Q_1$ such that $|Q_1\cap R|\ge |R|/r$. Let $R_1=R\setminus (Q_1\cap R)$. Then
$$|R_1|=|R|-|Q_1\cap R|\le |R|-\frac{|R|}{r}=\left(1-\frac{1}{r}\right)|R|,$$
and by Claim \ref{F/r}, there exists a perfect matching $Q_2$ such that $|Q_2\cap R_1|\ge |R_1|/r$. Let $R_2=R_1\setminus (Q_2\cap R_1)$. Then
$$|R_2|=|R_1|-|Q_2\cap R_1|\le |R_1|-\frac{|R_1|}{r}=\left(1-\frac{1}{r}\right)|R_1|\le \left(1-\frac{1}{r}\right)^2|R|.$$
Repeating this process $r$ times, we obtain
$$|R_r|\le \left(1-\frac{1}{r}\right)^r|R|\le \left(1-\frac{1}{r}\right)^r(8r+2),$$
where the last inequality follows from \eqref{8r+2}. We now prove the following inequality.
\begin{claim}
For every $r\ge1$, we have
$$\left(1-\frac{1}{r}\right)^r(8r+2)<3r.$$
\end{claim}
\begin{proof}
The case $r=1$ is immediate. For $r\ge2$, we use the logarithmic inequalities $ \ln(1-x)\le-x-x^2/2\;\;(0<x<1)$ and $\ln(1+y)\le y\;(y>0)$. Substituting $x=1/r$ and $y=1/(4r)$, we obtain
\begin{equation*}
r\ln\left(1-\frac1r\right)+\ln\left(1+\frac1{4r}\right)\le -1-\frac1{2r}+\frac1{4r}
=-1-\frac1{4r}< -1.
\end{equation*}
Consequently,
\begin{equation*}
\left(1-\frac1r\right)^r(8r+2)=8r\left(1-\frac1r\right)^r\left(1+\frac1{4r}\right)<\frac{8r}{e}<3r,
\end{equation*}
where $e$ denotes Euler's number. This proves the claim.
\end{proof}

So $|R_r|\le 3r-1$. If $|R_r|\le r$, then $R_r$ can be covered by the union of at most $r$ perfect matchings, since $G$ is matching-covered. Thus we conclude that
$$\tau(G)\le (r+1)+r+r=3r+1\le 4r.$$

Therefore, we may assume that $|R_r|=r+t$, where $1\le t\le 2r-1$. By Claim \ref{F/r}, there exists a perfect matching
$Q_{r+1}$ such that
$$|Q_{r+1}\cap R_r|\ge\left\lceil\frac{|R_r|}{r}\right\rceil=\left\lceil\frac{r+t}{r}\right\rceil
\ge\left\lceil\frac{r+1}{r}\right\rceil=2.$$
Choose two distinct edges $e_{r+1},f_{r+1}\in Q_{r+1}\cap R_r$, and let $R_{r+1}=R_r\setminus\{e_{r+1},f_{r+1}\}$. Then $|R_{r+1}|=r+t-2$. Repeating this step $\lceil t/2\rceil$ times, we have $|R_{r+\lceil t/2\rceil}|\le r$. Indeed,
\[
\big|R_{r+\lceil t/2\rceil}\big|=|R_r|-2\lceil t/2\rceil=r+t-2\lceil t/2\rceil
=
\begin{cases}
r & \text{if } t \text{ is even},\\
r-1 & \text{if } t \text{ is odd}.
\end{cases}
\]

If $t$ is odd, then $R_{r+\lceil t/2\rceil}$ can be covered by the union of at most $r-1$ perfect matchings. Thus we conclude that
$$\tau(G)\le r+1+r+(t+1)/2+r-1\le4r,$$
where the last inequality follows from $t\le 2r-1$.

If $t$ is even, then $R_{r+\lceil t/2\rceil}$ can be covered by the union of at most $r$ perfect matchings. Since $2r\ge t+1$ and $t+1$ is odd,  we have $2r\ge t+2$. Thus we conclude that
$$\tau(G)\le r+1+r+t/2+r\le4r.$$
This completes the proof of Theorem \ref{main-2}.

\vspace{2mm}
\noindent\textbf{Declarations}

\vspace{2mm}
\noindent{\bf Conflict of interest.} The authors declare that they have no conflict of interest.

\end{sloppypar}
\end{document}